\documentclass{amsart}

\usepackage{latexsym}
\usepackage{amsfonts}
\usepackage{amssymb}
\usepackage{amsmath}

\newcommand{\R}{\mathbb{R}}
\newcommand{\C}{\mathbb{C}}
\newcommand{\N}{\mathbb{N}}
\newcommand{\ord}{\textrm{ord}}

\newtheorem{theorem}{Theorem}
\newtheorem{lemma}[theorem]{Lemma}

\begin{document}

\title{A new Clunie type theorem for difference polynomials}

\author{Risto Korhonen}
\address{Department of Mathematics and Statistics, P.O. Box 68,
FI-00014 University of Helsinki, Finland}
\email{risto.korhonen@helsinki.fi}



\begin{abstract}
It is shown that if $w(z)$ is a finite-order meromorphic solution of the
equation
    \begin{equation*}
    H(z,w) P(z,w)=Q(z,w),
    \end{equation*}
where $P(z,w)=P(z,w(z),w(z+c_1),\ldots,w(z+c_n))$,
$c_1,\ldots,c_n\in\C$, is a homogeneous difference polynomial with
meromorphic coefficients, and $H(z,w)=H(z,w(z))$ and
$Q(z,w)=Q(z,w(z))$ are polynomials in $w(z)$ with meromorphic
coefficients having no common factors such that
    \begin{equation*}
    \begin{split}
    &\max\{\deg_w(H),\deg_w(Q)-\deg_w(P)\}> \min\{\deg_w(P),\ord_0(Q)\}-\ord_0(P),
    \end{split}
    \end{equation*}
where $\ord_0(P)$ denotes the order of zero of
$P(z,x_0,x_1,\ldots,x_n)$ at $x_0=0$ with respect to the variable
$x_0$, then the Nevanlinna counting function $N(r,w)$ satisfies
$N(r,w)\not=S(r,w)$. This implies that $w(z)$ has a relatively
large number of poles. For a smaller class of equations a stronger
assertion $N(r,w)=T(r,w)+S(r,w)$ is obtained, which means that the
pole density of $w(z)$ is essentially as high as the growth of
$w(z)$ allows. As an application, a simple necessary and
sufficient condition is given in terms of the value distribution
pattern of the solution, which can be used as a tool in ruling out
the possible existence of special finite-order Riccati solutions
within a large class of difference equations containing several
known difference equations considered to be of Painlev\'e type.
\end{abstract}

\maketitle

\section{Introduction and main results}

According to Clunie's theorem \cite{clunie:62}, if a meromorphic
function $f$ satisfies the differential equation
    \begin{equation}\label{clunieeq}
    f^nP(z,f)=Q(z,f),
    \end{equation}
where $n\in\N$, and $P(z,f)$ and $Q(z,f)$ are differential
polynomials in $f$ with meromorphic coefficients such that $\deg_f
Q(z,f)\leq n$, then the Nevanlinna proximity function
$m(r,\cdot\,)$ satisfies
    \begin{equation}
    m(r,P(z,f))=O(\log r + \log T(r,f)+ \mathcal{T}(r))\nonumber
    \end{equation}
where $r$ approaches infinity outside of a set of finite linear
measure, and $\mathcal{T}(r)$ is the maximum of the Nevanlinna
characteristics of the coefficients of $P(z,f)$ and $Q(z,f)$.
Originally Clunie used his result to consider certain properties
of entire and meromorphic functions, and later on Clunie's theorem
and its subsequent generalizations, see, e.g.,
\cite{doeringer:82,lahirib:04,yangy:07}, have proven to be
valuable tools in the study of value distribution of meromorphic
solutions of Painlev\'e, and other non-linear differential
equations, see, e.g., \cite{laine:93,gromakls:02}.

Ablowitz, Halburd and Herbst suggested that the existence of
sufficiently many finite-order meromorphic solutions is a unique
characteristic of a Painlev\'e type difference equation
\cite{ablowitzhh:00}. In \cite{halburdk:07PLMS} it was shown that
the existence of one finite-order meromorphic solution is enough
to reduce the second order difference equation
    \begin{equation}\label{class}
    w(z+1)+w(z-1)=R(z,w),
    \end{equation}
where $R(z,w)$ is rational in $w(z)$ with meromorphic
coefficients, into a list of equations consisting only of
difference Painlev\'e equations and linear equations within the
class (\ref{class}), provided that the finite-order solution
$w(z)$ does not simultaneously satisfy a difference Riccati
equation
    $$
    w(z+1)=\frac{a_1(z)w(z)+a_0(z)}{b_1(z)w(z)+b_0(z)}
    $$
where the coefficients are meromorphic functions (small with
respect to $w$ in the sense of Nevanlinna theory) such that
$a_1b_0\not\equiv a_0b_1$. An essential part of the method used in
this classification is based on a local analysis of the behavior
of a meromorphic solution near its poles, which can only be
performed non-vacuously if there are sufficiently many poles to
begin with. High pole density of solutions can be verified by
applying a direct difference analogue of Clunie's theorem
\cite{halburdk:06JMAA} (concerning the equation $f^nP(z,f)=Q(z,f)$
where $P(z,f)$ and $Q(z,f)$ are difference polynomials, and $f$ is
of finite order).

The class (\ref{class}) contains many equations considered to be
of Painlev\'e type, including three alternate versions of
difference Painlev\'e I, and a difference Painlev\'e II.
Nevertheless, most of the difference Painlev\'e equations fall
outside of class (\ref{class}) \cite{grammaticosnr:99}. For
instance, a known discretization of the Painlev\'e III equation,
    \begin{equation}\label{dpiii}
    w(z+1)w(z-1)=\frac{\gamma(z+1)w(z)^2+\zeta(z)\lambda^zw(z)+\mu(z)\lambda^{2z}}{(w(z)-1)(w(z)-\gamma(z))},
    \end{equation}
where $\lambda\in\C$, $\gamma$ and $\zeta$ are periodic
meromorphic functions with period two, and $\mu$ is a period one
meromorphic function, not only lies outside of the class
(\ref{class}) but the Clunie difference analogue in
\cite{halburdk:06JMAA} is inapplicable for this equation. This
causes a difficulty in making sure that solutions have enough
poles so that the local analysis needed for the classification can
be performed. Theorem~\ref{laineyangthm} below by I.~Laine and
C.~C.~Yang is a generalization of
\cite[Theorem~3.1]{halburdk:06JMAA} and applicaple for the
equation (\ref{dpiii}).

At this point we pause briefly to introduce the notation used in this paper. Let $c_j\in\C$ for $j=1,\ldots,n$ and let $I$ be a finite set of multi-indexes
$\lambda=(\lambda_0,\ldots,\lambda_n)$. A difference polynomial of
a meromorphic function $w(z)$ is defined as
    \begin{equation}\label{Pzw}
    \begin{split}
    P(z,w)&=P(z,w(z),w(z+c_1),\ldots,w(z+c_n))\\ & =\sum_{\lambda\in I}
    a_\lambda(z)w(z)^{\lambda_0} w(z+c_1)^{\lambda_1} \cdots
    w(z+c_n)^{\lambda_n},
    \end{split}
    \end{equation}
where the coefficients $a_\lambda(z)$ are \textit{small} with respect to
$w(z)$ in the sense that $T(r,a_\lambda)=o(T(r,w))$ as $r$ tends
to infinity outside of an exceptional set $E$ of finite
logarithmic measure
    \begin{equation}\label{finlogmeas}
    \lim_{r\to\infty}\int_{E\cap[1,r)} \frac{dt}{t}<\infty.
    \end{equation}
From now on we use an abbreviated notation $\int_E dt/t<\infty$
instead of \eqref{finlogmeas} to denote finite logarithmic
measure. The notation $T(r,a_\lambda)=S(r,w)$ is also used to
indicate that the characteristic function of $a_\lambda(z)$ is
small with respect to the characteristic of $w(z)$. The
\textit{total degree} of $P(z,w)$ in $w(z)$ and the shifts of
$w(z)$ is denoted by $\deg_w(P)$, and the \textit{order of a zero}
of $P(z,x_0,x_1,\ldots,x_n)$, as a function of $x_0$ at $x_0=0$,
is denoted by $\ord_0(P)$. (For instance, if
$P(z,w)=w(z)^2+w(z+1)w(z)$, then $\deg_w(P)=2$ and $\ord_0(P)=1$).
Moreover, the \textit{weight} of a difference polynomial
\eqref{Pzw} is defined by
    $$
    \kappa(P)=\max_{\lambda\in I}\left\{\sum_{j=1}^n \lambda_j\right\},
    $$
where $\lambda=(\lambda_0,\ldots,\lambda_n)$, and the set $I$ is
the same as in (\ref{Pzw}) above. By this definition the weight of
the polynomial $w(z+1)w(z-1)+w(z)w(z-1)+w(z)w(z+1)$ is two, for
instance. The difference polynomial $P(z,w)$ is said to
be \textit{homogeneous} with respect to $w(z)$ if the degree
$d_\lambda=\lambda_0+\cdots+\lambda_n$ of each term in the sum
(\ref{Pzw}) is non-zero and the same for all $\lambda\in I$. Finally, the \textit{order of growth} of a meromorphic function $w$ is defined by
    \begin{equation*}
    \rho(w)=\limsup_{r\to\infty}\frac{\log T(r,w)}{\log r}.
    \end{equation*}
Notation and fundamental results from Nevanlinna theory are
frequently used throughout this paper, see, e.g.,
\cite{cherryy:01,goldbergo:08,hayman:64}.

\begin{theorem}[Laine, Yang \cite{lainey:07}]\label{laineyangthm}
Let $f$ be a transcendental meromorphic solution of finite order
$\rho$ of a difference equation of the form
    \begin{equation}
    U(z, f)P(z, f) = Q(z, f),\nonumber
    \end{equation}
where $U(z, f)$,  $P(z, f)$ and $Q(z, f)$ are difference
polynomials such that the total degree $\deg_f U(z, f) = n$ in
$f(z)$ and its shifts, and $\deg_fQ(z, f) \leq n$. If $U(z, f)$
contains just one term of maximal total degree in $f(z)$ and its
shifts, then, for each $\varepsilon>0$,
    \begin{equation}
    m(r, P(z, f)) = O(r^{\rho-1+\varepsilon}) + S(r, f),\nonumber
    \end{equation}
possibly outside of an exceptional set of finite logarithmic
measure.
\end{theorem}

With the help of Theorem~\ref{laineyangthm} the full
classification of
   \begin{equation}\label{class2}
    w(z+1)w(z-1)=R(z,w)
    \end{equation}
containing (\ref{dpiii}) has been completed in
\cite{ronkainen:09}. Although Theorem~\ref{laineyangthm} covers a
large class of equations, the equation known as the difference
Painlev\'e IV (d-P$_{IV}$) is not one of them (see
\cite{ramanigh:91} for a discretization of the Painlev\'e IV
equation). One of the aims of this paper is to prove the following
alternative version of the difference Clunie lemma for a class of
difference equations which includes d-P$_{IV}$.

\begin{theorem}\label{1stthm}
Let $w(z)$ be a finite-order meromorphic solution of
    \begin{equation}\label{cleq0}
    H(z,w) P(z,w)=Q(z,w),
    \end{equation}
where $P(z,w)$ is a homogeneous difference polynomial with
meromorphic coefficients, and $H(z,w)$ and $Q(z,w)$ are
polynomials in $w(z)$ with meromorphic coefficients having no
common factors. If
    \begin{equation}\label{cleqassumpt}
    \begin{split}
    &\max\{\deg_w(H),\deg_w(Q)-\deg_w(P)\}> \min\{\deg_w(P),\emph{\ord}_0(Q)\}-\emph{\ord}_0(P),
    \end{split}
    \end{equation}
then $N(r,w)\not=S(r,w)$.
\end{theorem}

The expression $N(r,w)\not=S(r,w)$ means that the pole counting
function of $w$ is not small with respect to the characteristic
function of $w$. In other words, there exists an absolute constant
$K\in (0,1]$ and a set $E$ of infinite logarithmic measure, such
that $N(r,w)\geq K\,T(r,w)$ for all $r\in E$.

The class of equations (\ref{cleq0}) contains many difference
equations considered to be of Painlev\'e type, including equations
known as difference Painlev\'e I--IV, for suitable choices of the
polynomials $H$, $P$ and $Q$. In Section \ref{sec2} below we will
consider a class of equations within (\ref{cleq0}) containing
d-P$_{IV}$ as an example of an application of
Theorem~\ref{1stthm}.

By adding a constraint to the degrees and weights of the
difference polynomials in Theorem~\ref{1stthm}, the following
stronger assertion is obtained.

\begin{theorem}\label{2ndthm}
Let $w(z)$ be a finite-order meromorphic solution of
    \begin{equation}\label{cleq2}
    H(z,w) P(z,w)=Q(z,w),
    \end{equation}
where $P(z,w)$ is a homogeneous difference polynomial with
meromorphic coefficients, and $H(z,w)$ and $Q(z,w)$ are
polynomials in $w(z)$ with meromorphic coefficients having no
common factors. If
    \begin{equation}\label{condcleq2}
    2\kappa(P) \leq \max\{\deg_w(Q),\deg_w(H)+\deg_w(P)\}-\min\{\deg_w(P),\emph{\ord}_0(Q)\},
    \end{equation}
then, for any $\delta\in(0,1)$,
    \begin{equation}\label{mvsT}
    m(r,w)=o\left(\frac{T(r,w)}{r^\delta}\right) + O(\mathcal{T}(r)),
    \end{equation}
where $r$ runs to infinity outside of an exceptional set of finite
logarithmic measure, and $\mathcal{T}(r)$ is the maximum of the
Nevanlinna characteristics of the coefficients of $P(z,w)$,
$Q(z,w)$ and $H(z,w)$.
\end{theorem}

Theorems~\ref{1stthm} and \ref{2ndthm} provide some of the
necessary tools needed to single out Painlev\'e type equations out
of the difference equation
    \begin{equation}\label{exampleeqIV}
    \overline{w}\underline{w}+\overline{w}w+w\underline{w}=\frac{a_3w^3+a_2w^2+a_1w+a_0}{(w-b)(w-c)}
    \end{equation}
where the coefficients are rational functions, $a_0\not\equiv0$,
and we have suppressed the $z$-dependence of $w(z)$ by writing
$\overline{w}\equiv w(z+1)$, $\underline{w}\equiv w(z-1)$ and
$w\equiv w(z)$. Namely, by taking $H(z,w)=(w-b)(w-c)$,
$Q(z,w)=a_3w^3+a_2w^2+a_1w+a_0$ and
$P(z,w)=\overline{w}\underline{w}+\overline{w}w+w\underline{w}$ in
Theorem~\ref{2ndthm}, it follows that $\kappa(P)=2$,
$\deg_w(Q)=3$, $\deg_w(P)=2$, $\deg_w(H)=2$ and $\ord_0(Q)=0$.
Hence, the assumption \eqref{condcleq2} is satisfied, and so
$m(r,w)=S(r,w)$ by Theorem~\ref{2ndthm}. Therefore,
$N(r,w)=T(r,w)+S(r,w)$, and all non-rational finite-order
meromorphic solutions of \eqref{exampleeqIV} have nearly as many
poles as their growth enables. In particular, this is true for all
non-rational finite-order meromorphic solutions of d-P$_{IV}$,
since this equation is a special case of \eqref{exampleeqIV}.


\section{Applications to difference Riccati equation}\label{sec2}

The Painlev\'e property has proved to be a good detector of
integrability in differential equations \cite{ablowitzc:91}. In
the beginning of the $20^\textrm{th}$ century, Painlev\'e, Gambier
and Fuchs identified all those equations that possess this
property out of a large class of second-order ordinary
differential equations
\cite{fuchs:05,gambier:10,painleve:00,painleve:02}. All of the
equations could be solved in terms of previously known functions,
solutions of linear equations, or in terms of solutions of one of
six new equations, now known as the Painlev\'e equations. The
Painlev\'e equations were later on integrated by using inverse
scattering transform techniques, see, e.g.,~\cite{ablowitzs:77}.
In the first order case Malmquist \cite{malmquist:13} has shown
that the existence of one meromorphic solution of the differential
equation
    \begin{equation}\label{malmq}
    w'=R(z,w),
    \end{equation}
where $R(z,w)$ is rational in both arguments, reduces
\eqref{malmq} into a Riccati equation. A simple proof of this fact
was given later on by Yosida \cite{yosida:33} using techniques
from Nevanlinna theory.

There are several candidates for the discrete Painlev\'e property,
including the singularity confinement by Grammaticos, Ramani and
Papageorgiou \cite{grammaticosrp:91}, algebraic entropy by
Hietarinta and Viallet \cite{hietarintav:98}, the existence of
sufficiently many finite-order meromorphic solutions by Ablowitz,
Halburd and Herbst \cite{ablowitzhh:00}, and Diophantine
integrability by Halburd~\cite{halburd:05}. As mentioned in the
introduction, Halburd and the author \cite{halburdk:07PLMS} showed
that the existence of one finite-order meromorphic solution
growing faster than the coefficients is sufficient to reduce the
second order difference equation (\ref{class}), where $R(z,w)$ is
rational in $w(z)$ with meromorphic coefficients, into a list of
equations consisting of difference Painlev\'e equations and linear
equations within the class (\ref{class}), provided that the
finite-order solution $w(z)$ does not simultaneously satisfy a
difference Riccati equation. Most of the equations of Painlev\'e
type - both continuous and discrete - are known to possess special
solutions satisfying a first order Riccati equation for particular
choices of their parameter values
\cite{gromakls:02,tamizhmanitgr:04}. Although this property is
considered to be one of the typical characteristics of Painlev\'e
equations, Riccati type solutions also appear as special solutions
of non-integrable equations. Therefore, in the
Ablowitz-Halburd-Herbst approach, the existence of any number of
finite-order meromorphic solutions is insufficient to indicate
integrability of a difference equation, if these solutions happen
to be simultaneously solutions to a Riccati equation. At the same
time, the existence of already one non-Riccati finite-order
solution appears to indicate a Painlev\'e type difference
equation.

The purpose of this section is to give a simple necessary and
sufficient condition which can be used to rule out the possible
existence of special finite-order Riccati type solutions within a
large class of difference equations. The condition is formulated
in terms of the value distribution pattern of the considered
meromorphic solution near its poles, which is straightforward to
work out for most difference equations.

\begin{theorem}\label{mainthm}
Let $a(z)$ and $c(z)$ be rational functions, and let $w(z)$ be a
non-rational finite-order meromorphic solution of
    \begin{equation}\label{cleq}
    H(z,w) P(z,w)=Q(z,w),
    \end{equation}
where $P(z,w)$ is a homogeneous difference polynomial with respect
to $w(z)$ having rational coefficients, and $H(z,w)$ and $Q(z,w)$
are polynomials in $w(z)$ with rational coefficients having no
common factors. If
    \begin{equation}\label{cond2}
    2\kappa(P) \leq \max\{\deg_w(Q),\deg_w(H)+\deg_w(P)\}-\min\{\deg_w(P),\emph{\ord}_0(Q)\},
    \end{equation}
then the following statements are equivalent:
    \begin{itemize}
    \item[(i)] There exists a
    positive integer $k_{\hat{z}}$ and complex constants
    $\alpha_{\hat{z}}$, $\beta_{\hat{z}}\not=0$ and $\gamma_{\hat{z}}$ such that, at
    all except at most finitely many poles $\hat{z}$ of $w(z)$,
        \begin{eqnarray*}
        w(z-1)&=&c(z-1)+\alpha_{\hat{z}}(z-\hat{z})^{k_{\hat{z}}}+O((z-\hat{z})^{k_{\hat{z}}+1})\\
        w(z)&=&\beta_{\hat{z}}(z-\hat{z})^{-k_{\hat{z}}}+O((z-\hat{z})^{1-k_{\hat{z}}})\\
        w(z+1)&=&a(z)+\gamma_{\hat{z}}(z-\hat{z})^{k_{\hat{z}}}+O((z-\hat{z})^{k_{\hat{z}}+1})
        \end{eqnarray*}
    for all $z$ in a neighborhood of $\hat{z}$.
    \item[(ii)] The function $w(z)$ is a solution of the difference Riccati equation
        \begin{equation}\label{Riccati0}
        w(z+1)=\frac{a(z)w(z)+b(z)}{w(z)-c(z)},
        \end{equation}
    where $b(z)$ is a meromorphic function having at most finitely
    many poles, and satisfying $\rho(b)\leq\max\{0,\rho(w)-1\}$.
    \end{itemize}
\end{theorem}

To demonstrate the use of Theorem~\ref{mainthm}, we consider the
(in general non-integrable) non-autonomous difference equation
    \begin{equation}\label{exampleeq0}
    \overline{w}+\underline{w}=\frac{a_2 w^2+a_1 w +a_0}{w^2}
    \end{equation}
which contains a known discrete form of the Painlev\'e I (see,
e.g., \cite{fokasgr:93}) as a special case. Assuming that $w$ is a
non-rational finite-order meromorphic solution of
(\ref{exampleeq0}) where the coefficients are rational functions
and $a_0\not\equiv0$, it follows by applying Theorem~\ref{2ndthm}
with $P(z,w)=\overline{w}+\underline{w}$, $H(z,w)=w^2$ and
$Q(z,w)=a_2w^2+a_1w+a_0$ that $m(r,w)=S(r,w)$. Suppose that $w$
has only finitely many poles. Then, since $N(r,w)=O(\log r)$ and
$T(r,w)=N(r,w)+S(r,w)$, it follows that $T(r,w)=S(r,w)+O(\log r)$.
But this is impossible due to the fact that $w$ is non-rational,
and so $w$ has infinitely many poles. Suppose now that $w$ has
finitely many zeros. Then $N(r,1/w)=O(\log r)$, and, by
Lemma~\ref{logdiff} below,
    \begin{equation*}
    \begin{split}
    m(r,\overline{w}+\underline{w})&=m\left(r,\frac{w(\overline{w}+\underline{w})}{w}\right)
    \leq m(r,w)+m\left(r,\frac{\overline{w}+\underline{w}}{w}\right)\\ &= m(r,w)+S(r,w)=S(r,w).
    \end{split}
    \end{equation*}
Therefore, by the Valiron-Mohon'ko identity
\cite{valiron:31,mohonko:71} (see also
\cite[Theorem~6.5]{goldbergo:08} and
\cite[Theorem~2.2.5]{laine:93}), it follows that
    \begin{equation*}
    \begin{split}
    2T(r,w)&=T(r,\overline{w}+\underline{w})+O(\log r)\\
    &=m(r,\overline{w}+\underline{w})+N(r,\overline{w}+\underline{w}) +O(\log r)\\
    &=N(r,\overline{w}+\underline{w})+S(r,w)+O(\log r).\\
    \end{split}
    \end{equation*}
Since by \eqref{exampleeq0},
    \begin{equation*}
    N(r,\overline{w}+\underline{w})=N\left(r,\frac{1}{w^2}\right)+O(\log r)
    =2N\left(r,\frac{1}{w}\right)+O(\log r)=O(\log r),
    \end{equation*}
it follows that $T(r,w)=S(r,w)+O(\log r)$ which is impossible.
Therefore $w$ has also infinitely many zeros. By equation
(\ref{exampleeq0}) it follows that whenever $w$ has a zero of
multiplicity $k_0$ at $z=z_0$, then $w$ has a pole at least of
multiplicity $2k_0$ at $z=z_0+1$ or $z=z_0-1$. By
Theorem~\ref{mainthm} these poles are not of the type allowed for
the solution of a Riccati difference equation. Since we have shown
that $w$ has infinitely many such poles, it follows by
Theorem~\ref{mainthm} that $w$ cannot be a solution of the Riccati
difference equation (\ref{Riccati0}).

As another example we consider the equation
    \begin{equation}\label{exampleeq}
    \overline{w}\underline{w}+\overline{w}w+w\underline{w}=\frac{a_3w^3+a_2w^2+a_1w+a_0}{(w-b)(w-c)}
    \end{equation}
where the coefficients are rational functions, and
$a_0\not\equiv0$. Equation (\ref{exampleeq}) contains a known
discrete form of the Painlev\'e IV (see, e.g., \cite{ramanigh:91})
as a special case, which is known to have special difference
Riccati solutions. If $w$ is a non-rational finite-order solution
of (\ref{exampleeq}), then, similarly as above, it follows by
applying Theorem~\ref{2ndthm} with
$P(z,w)=\overline{w}\underline{w}+\overline{w}w+w\underline{w}$,
$H(z,w)=(w-b)(w-c)$ and $Q(z,w)=a_3w^3+a_2w^2+a_1w+a_0$ that there
are infinitely many points $\hat{z}$ where $w(\hat{z})=\infty$.
By substituting a suitable Laurent series expansion into (\ref{exampleeq}) it follows that all except
possibly finitely many of these poles appear as a part of one of
the following sequences:
    \begin{eqnarray}
    && w(\hat{z}-1)=b(\hat{z}-1),\quad w(\hat{z})=\infty,  \quad w(\hat{z}+1)=a_3(\hat{z})-b(\hat{z}-1) \label{seq1}\\
    && w(\tilde{z}-1)=c(\tilde{z}-1),\quad w(\tilde{z})=\infty,  \quad w(\tilde{z}+1)=a_3(\tilde{z})-c(\tilde{z}-1) \label{seq2}\\
    && w(\breve{z}-1)=a_3(\breve{z})-b(\breve{z}+1),\quad w(\breve{z})=\infty,  \quad w(\breve{z}+1)=b(\breve{z}+1) \label{seq3}\\
    && w(\acute{z}-1)=a_3(\acute{z})-c(\acute{z}+1),\quad w(\acute{z})=\infty,  \quad w(\acute{z}+1)=c(\acute{z}+1) \label{seq4}\\
    && w(\grave{z}-1)=\infty,\quad w(\grave{z})=K_{\grave{z}},\quad w(\grave{z}+1)=\infty
    \label{seq5},
    \end{eqnarray}
where $\hat{z}, \tilde{z}, \breve{z}, \acute{z}, \grave{z} \in \C$
and $K_{\grave{z}}\in\C \cup \{\infty\}$. (The possible finitely
many exceptional poles which are not one of the types
(\ref{seq1})--(\ref{seq4}) arise from the poles of the
coefficients.) According to Theorem~\ref{mainthm}, the function
$w(z)$ is a special Riccati solution of (\ref{exampleeq}) if and
only if all except possibly finitely many poles of $w(z)$ are of
exactly one of the types (\ref{seq1})--(\ref{seq4}).

\section{Proofs of theorems}\label{proof1}

We begin by stating a known difference analogue of the lemma on
the logarithmic derivative.

\begin{lemma}[\cite{halburdk:06JMAA,halburdk:06AASFM}]\label{logdiff}
Let $f$ be a non-constant meromorphic function of finite order,
$c\in\C$ and $\delta\in(0,1)$. Then
    \begin{equation}\label{diff}
    m\left(r,\frac{f(z+c)}{f(z)}\right)= o\left(\frac{T(r,f)}{r^\delta}\right)
    \end{equation}
where $r$ approaches infinity outside of a possible exceptional
set $E$ with finite logarithmic measure
$\int_E\frac{dr}{r}<\infty$.
\end{lemma}

See \cite[Corollary 2.5]{chiangf:08} for an alternative version of
Lemma \ref{logdiff}. The following generalization of \cite[Lemma
2.1]{halburdk:07PLMS} is needed in the proof Theorem~\ref{2ndthm}.

\begin{lemma}\label{technical}
Let $T:[0,+\infty)\to[0,+\infty)$ be a non-decreasing continuous
function, let $\delta\in(0,1)$, and let $s\in(0,\infty)$. If $T$
is of finite order, i.e.,
    \begin{equation}\label{assu}
    \limsup_{r\to\infty}\frac{\log T(r)}{\log r}<\infty,
    \end{equation}
then
   \begin{equation}
    T(r+s) = T(r)+ o\left(\frac{T(r)}{r^\delta}\right)\nonumber
    \end{equation}
where $r$ runs to infinity outside of a set of finite logarithmic
measure.
\end{lemma}

Note that by using \cite[Theorem 2.2 and Corollary
2.5]{chiangf:08} instead of Lemma~\ref{logdiff} the proof of
Theorem~\ref{2ndthm} could be simplified in the sense that
Lemma~\ref{technical} would no longer be required. However, this
would change the assertion \eqref{mvsT} of Theorem~\ref{2ndthm}
into
    \begin{equation}\label{cf}
    m(r,w)=O(r^{\rho(w)-1+\varepsilon}) + O(\mathcal{T}(r)),
    \end{equation}
where $\rho(w)$ is the order of $w$, $\varepsilon>0$ and $r>0$.
Unfortunately \eqref{cf} does not necessarily imply that $m(r,w)$
is small compared to $T(r,w)$ for all, or even most values of $r$.
Namely, if the \textit{lower order} of $w$, defined by
    $$
    \mu(w)=\liminf_{r\to\infty}\frac{\log T(r,w)}{\log r}
    $$
satisfies $\mu(w)<\rho(w)-1$, then $T(r,w)<r^{\rho(w)-1}$ in a
significant (and in some cases the largest) part of
$\mathbb{R}_+$. Therefore, for these particular values of $r$,
equation~\eqref{cf} gives no information on the relative size of
$m(r,w)$ compared to $T(r,w)$, and so the set where $m(r,w)$ may
not be small compared to $T(r,w)$ in \eqref{cf} can be much larger
than the exceptional set in \eqref{mvsT}. On the other hand, if
the growth of $w$ is assumed to be sufficiently regular in the
sense that $\mu(w)>\rho(w)-1$, then \eqref{cf} implies that
$m(r,w)=o(T(r,w)) + O(\mathcal{T}(r))$ without an exceptional set.

\medskip

\textit{Proof of Lemma \ref{technical}: } Denote
$\nu(r)=1-r^{-\delta}$, and assume conversely to the assertion
that the set $F\subset\R^{+}$ of all $r$ such that
    \begin{equation}
    T(r) \leq \nu(r) T(r+s)\nonumber
    \end{equation}
is of infinite logarithmic measure. Set $r_n=\min\{F\cap
[r_{n-1}+s,\infty)\}$ for all $n\in\N$, where $r_0$ is the
smallest element of $F$. Then the sequence $\{r_n\}_{n\in\N}$
satisfies $r_{n+1}-r_n\geq s$ for all $n\in\N$, $F\subset
\bigcup_{n=0}^\infty [r_n,r_n+s]$ and
    \begin{equation}\label{assuinpr}
    T(r_n) \leq \nu(r_n) T(r_{n+1})
    \end{equation}
for all $n\in\N$.

Let $\varepsilon\in(0,\delta^{-1}-1)$, and suppose that there
exist an $m\in\N$ such that $r_n\geq n^{1+\varepsilon}$ for all
$r_n\geq m$. But then,
    \begin{eqnarray*}
    \int_F\frac{dt}{t} &\leq& \sum_{n=0}^\infty
    \int_{r_n}^{r_n+s}\frac{dt}{t}
    \leq \int_1^{m} \frac{dt}{t} +  \sum_{n=1}^\infty
    \log\left(1+\frac{s}{r_n}\right)\\
    &\leq& \sum_{n=1}^\infty
    \log\left(1+s n^{-(1+\varepsilon)}\right) +O(1)  <\infty
    \end{eqnarray*}
which contradicts the assumption $\int_F\frac{dt}{t}=\infty$.
Therefore the sequence $\{r_n\}_{n\in\N}$ has a subsequence
$\{r_{n_j}\}_{j\in\N}$ such that $r_{n_j}\leq n_j^{1+\varepsilon}$
for all $j\in\N$. By iterating~(\ref{assuinpr}) along the sequence
$\{r_{n_j}\}$ and using the fact that $\nu(r)$ is an increasing
function, it follows that
    \begin{equation}
    T(r_{n_j}) \geq \frac{1}{\nu(r_j)^{n_j}} T(r_0)\nonumber
    \end{equation}
for all $j\in \N$, and hence
    \begin{eqnarray*}
    \limsup_{r\to\infty}\frac{\log T(r)}{\log r}&\geq& \limsup_{j\rightarrow\infty}
    \frac{\log T(r_{n_j})}{\log r_{n_j}}\\
    &\geq & \limsup_{j\rightarrow\infty} \frac{-n_j\log \nu(r_j)+\log T(r_0)}
    {\log r_{n_j}}\\
    &\geq &\limsup_{j\rightarrow\infty} \frac{-n_j\log (1-n_j^{-(1+\varepsilon)\delta})+\log T(r_0)}
    {(1+\varepsilon)\log n_j}=\infty
    \end{eqnarray*}
since $(1+\varepsilon)\delta<1$. This contradicts (\ref{assu}),
and so the assertion follows. \hfill $\Box$\par\vspace{2.5mm}

\textit{Proof of Theorem \ref{1stthm}: } Taking into account the fact that $P(z,w)$
is homogeneous, it follows by Lemma \ref{logdiff} that
    \begin{equation}\label{th1eq}
    m\left(r,\frac{P(z,w)}{w^{\deg_w(P)}}\right)=o\left(\frac{T(r,f)}{r^\delta}\right)
    \end{equation}
for any $\delta\in(0,1)$, and for all $r$ outside of an
exceptional set of finite logarithmic measure. Moreover, by
applying an identity due to Valiron \cite{valiron:31} and Mohon'ko
\cite{mohonko:71} (see also \cite[Theorem~6.5 and Appendix~B,
p.~453]{goldbergo:08}) to (\ref{cleq0}), it follows that
    \begin{equation}
    \begin{split}
    T\left(r,\frac{P(z,w)}{w^{\deg_w(P)}}\right)
    &= d_w T(r,w) + O(\mathcal{T}(r)),
    \label{T}
    \end{split}
    \end{equation}
where
    \begin{equation}\label{dw}
    d_w=\max\{\deg_w(Q),\deg_w(H)+\deg_w(P)\}-\min\{\deg_w(P),\ord_0(Q)\}
    \end{equation}
and $r$ approaches infinity outside of an exceptional set of finite logarithmic measure.
By combining (\ref{th1eq}), (\ref{T}) and (\ref{cleqassumpt}) it follows that
    \begin{equation}\label{th1eq2}
    N\left(r,\frac{P(z,w)}{w^{\deg_w(P)}}\right)\geq (1+\deg_w(P)-\ord_0(P))T(r,w)+S(r,w).
    \end{equation}

Suppose now on the contrary to the assertion of
Theorem~\ref{1stthm} that $N(r,w)=S(r,w)$. Therefore, denoting
$C=\max_{j=1,\ldots,n}\{|c_j|\}$ in (\ref{Pzw}), it follows from
Lemma \ref{technical} that
    \begin{equation*}
    \begin{split}
    N\left(r,\frac{P(z,w)}{w^{\ord_0(P)}}\right)&\leq (\deg_w(P)-\ord_0(P)) N(r+C,w)\\
    &=(\deg_w(P)-\ord_0(P)) N(r,w)+S(r,w)\\&=S(r,w).
    \end{split}
    \end{equation*}
Thus,
    \begin{equation}
    \begin{split}
    N\left(r,\frac{P(z,w)}{w^{\deg_w(P)}}\right)&\leq  N\left(r,\frac{P(z,w)}{w^{\ord_0(P)}}\right)
    +N\left(r,\frac{1}{w^{\deg_w(P)-\ord_0(P)}}\right)\\
    &= N\left(r,\frac{1}{w^{\deg_w(P)-\ord_0(P)}}\right) + S(r,w)\\
    &\leq T\left(r,\frac{1}{w^{\deg_w(P)-\ord_0(P)}}\right)+S(r,w).\nonumber
    \end{split}
    \end{equation}
Hence, by the first main theorem of Nevanlinna theory, it follows
that
    \begin{equation}
    N\left(r,\frac{P(z,w)}{w^{\deg_w(P)}}\right)\leq (\deg_w(P)-\ord_0(P))T(r,w)+S(r,w)\nonumber
    \end{equation}
which contradicts (\ref{th1eq2}). We conclude that
$N(r,w)\not=S(r,w)$. \hfill $\Box$\par\vspace{2.5mm}

\textit{Proof of Theorem \ref{2ndthm}: } Suppose now that $w(z)$
is a finite-order meromorphic solution of (\ref{cleq2}) such that
(\ref{condcleq2}) holds. By denoting
$C=\max_{j=1,\ldots,n}\{|c_j|\}$ in (\ref{Pzw}), it follows that
    \begin{equation}\label{aux1}
    N\left(r,\frac{P(z,w)}{w^{\deg_w(P)}}\right)\leq
    \kappa(P)\left(N(r+C,w)+N\left(r,\frac{1}{w}\right)\right)+O(\mathcal{T}(r)).
    \end{equation}
Since by Lemma \ref{technical}
    \begin{equation}
    N(r+C,w)=N(r,w)+o\left(\frac{N(r,w)}{r^\delta}\right)\nonumber
    \end{equation}
for all $r$ outside of a $E$ set of finite logarithmic measure,
inequality (\ref{aux1}) yields
    \begin{equation}\label{aux2}
    N\left(r,\frac{P(z,w)}{w^{\deg_w(P)}}\right)\leq \kappa(P)(2T(r,w)-m(r,w))+o\left(\frac{T(r,w)}{r^\delta}\right) +
    O(\mathcal{T}(r))
    \end{equation}
for all $r\not\in E$. On the other hand, by (\ref{th1eq}) and
(\ref{T}),
    \begin{equation}\label{aux3}
    \begin{split}
    N\left(r,\frac{P(z,w)}{w^{\deg_w(P)}}\right)= d_w T(r,w)
     +o\left(\frac{T(r,w)}{r^\delta}\right) + O(\mathcal{T}(r)),
    \end{split}
    \end{equation}
where $r$ lies outside of a set $F$ of finite logarithmic measure,
and $d_w$ is as in \eqref{dw}. By combining inequalities
(\ref{aux2}) and (\ref{aux3}) with the assumption
(\ref{condcleq2}), it follows that
    \begin{equation*}
    m(r,w)=o\left(\frac{T(r,w)}{r^\delta}\right) + O(\mathcal{T}(r))
    \end{equation*}
for all $r\not\in E\cup F$. \hfill $\Box$\par\vspace{2.5mm}

\textit{Proof of Theorem \ref{mainthm}: } Suppose first that all
except finitely many poles of $w(z)$ are in a sequence of the type
(i). Recall that we have adopted the short notation $w=w(z)$,
$\overline{w}=w(z+1)$ and $\underline{w}=w(z-1)$. The auxiliary
function $g$ defined by
    \begin{equation}\label{g}
    g=(\overline{w}-a)(w-c)
    \end{equation}
is meromorphic, of finite order, and, by Lemma \ref{logdiff} and
Theorem~\ref{2ndthm}, it satisfies
    \begin{equation}
    \begin{split}
    m(r,g)&\leq m(r,\overline{w})+m(r,w)+m(r,a)+m(r,c)+O(1) \\
    &\leq 2m(r,w)+m\left(r,\frac{\overline{w}}{w}\right) +O(\log r)
    \\
    &= o\left(\frac{T(r,w)}{r^\delta}\right)+O(\log r)
    \label{mrg}
    \end{split}
    \end{equation}
as $r\to\infty$ outside of a set of finite logarithmic measure.
Moreover, all possible poles of $g$, with at most finitely many
exceptions, arise from poles of $w$ or $\overline{w}$ and
therefore are part of the sequence in (i). Suppose first that
$w(\hat{z})=\infty$ with multiplicity $k_{\hat{z}}$. Then by the
sequence in (i), $w(\hat{z}+1)=a(\hat{z})$ with multiplicity no
less than $k_{\hat{z}}$, and so $g$ assumes a finite value at
$z=\hat{z}$. Similarly, if $w(\hat{z}+1)=\infty$ with multiplicity
$k_{\hat{z}+1}$, then $w(\hat{z})=c(\hat{z})$ with multiplicity no
less than $k_{\hat{z}+1}$, and so $g$ is again finite at
$z=\hat{z}$. Hence $g$ has only finitely many poles. By combining
this fact with (\ref{mrg}), it follows that
    \begin{equation}
    T(r,g)=o\left(\frac{T(r,w)}{r^\delta}\right)+O(\log r)\nonumber
    \end{equation}
for all $r$ outside of an exceptional set of finite logarithmic
measure. Therefore,
    \begin{equation}\label{estim}
    T(r,g)\leq r^{\rho(w)+\varepsilon-\delta}+K\log r
    \end{equation}
where $K>0$ is an absolute constant, $\varepsilon>0$, and $r$ lies
outside an exceptional set of finite logarithmic measure. By
\cite[Lemma 1.1.2]{laine:93} the exceptional set can be removed if
$r$ is replaced by $r^{1+\varepsilon}$ on the right side of
(\ref{estim}). Then (\ref{estim}) becomes
    \begin{equation}
    T(r,g)\leq r^{(1+\varepsilon)(\rho(w)+\varepsilon-\delta)}+K(1+\varepsilon)\log r\nonumber
    \end{equation}
for all $r$ sufficiently large. Since $\varepsilon>0$ and
$\delta\in(0,1)$ are arbitrary, it follows that
$\rho(g)\leq\max\{0,\rho(w)-1\}$. The first part of the
assertion follows by choosing $b=g-ac$.

Assume now that $w(z)$ is a finite-order meromorphic solution of
the Riccati equation \eqref{Riccati0} in (ii), and suppose that $w(z)$ has a pole
of order $k_{\hat{z}}$ at $z=\hat{z}$. Then, by writing the Riccati equation \eqref{Riccati0} in the form
    \begin{equation*}
    w(z+1)=a(z)+\frac{a(z)c(z)+b(z)}{w(z)-c(z)}
    \end{equation*}
and substituting $z=\hat{z}$, it follows that
either
    \begin{equation}
     w(z+1)=a(z)+\gamma_{\hat{z}}(z-\hat{z})^{k_{\hat{z}}}+O((z-\hat{z})^{k_{\hat{z}}+1})
     \nonumber
    \end{equation}
for all $z$ in a small enough neighborhood of $\hat{z}$, or
$a(z)c(z)+b(z)$ and/or $c(z)$ has a pole at $z=\hat{z}$. The
former case is, as required, one of the entries of the sequence in
(i), while the latter case can occur only at most finitely many
times.

By writing the difference Riccati equation in (ii) as
    \begin{equation}
    w(z-1)=\frac{c(z-1)w(z)+b(z-1)}{w(z)-a(z-1)}\nonumber
    \end{equation}
it follows, similarly as above, that either
    \begin{equation}\label{finalentry}
     w(z-1)=c(z-1)+\alpha_{\hat{z}}(z-\hat{z})^{k_{\hat{z}}}+O((z-\hat{z})^{k_{\hat{z}}+1})
    \end{equation}
for all $z$ in a small enough neighborhood of $\hat{z}$, or $a(z)c(z)+b(z)$
and/or $a(z)$ has a pole at $z=\hat{z}-1$. As above we conclude
that (\ref{finalentry}) holds for all except finitely many poles
of $w$. \hfill $\Box$\par\vspace{2.5mm}

\section*{Acknowledgements}

The research reported in this paper was supported in part by the
Academy of Finland grant \#118314 and \#210245. We would like to
thank the anonymous referees for their helpful comments on the paper.


\def\cprime{$'$}

\end{document}